\documentclass[twoside, a4paper, 10pt, reqno]{amsart}
\usepackage{amscd}      
\usepackage{amssymb} 
\usepackage{amsfonts}

\newtheorem{theorem}{Theorem}[section]
\newtheorem{lemma}[theorem]{Lemma}
\newtheorem{proposition}[theorem]{Proposition}
 
\theoremstyle{definition}  
\newtheorem{definition}[theorem]{Definition}

\newtheorem{remark}[theorem]{Remark}

\newcommand{\DGQUE}{\operatorname{DGQUE}} 
\newcommand{\Ob}{\operatorname{Ob}} 
\newcommand{\DGLBA}{\operatorname{DGLBA}} 
\newcommand{\pr}{\operatorname{pr}} 
\newcommand{\id}{\operatorname{id}} 
\newcommand{\inc}{\operatorname{inc}} 
\newcommand{\Ker}{\text{Ker\,}} 
\newcommand{\Alt}{\operatorname{Alt}}

\newcommand{\CYB}{\operatorname{CYB}}

\newcommand{\G}{\operatorname{G}} 
\newcommand{\fd}{\operatorname{fd}} 
\newcommand{\inv}{\operatorname{inv}} 
\newcommand{\poly}{\operatorname{poly}} 
\newcommand{\cH}{\operatorname{cH}} 
\newcommand{\Lie}{\operatorname{Lie}} 
\newcommand{\alg}{\operatorname{alg}} 
\newcommand{\Gi}{\G_\infty} 
\newcommand{\Ger}{\operatorname{Ger}}

\newcommand{\DQ}{\operatorname{DQ}}
\newcommand{\LGG}{\operatorname{L-G}}
\newcommand{\LG}{\LGG-G_\infty}
\newcommand{\LGi}{\LGG_\infty}

\newcommand{\op}{{\text{op}}}

\newcommand{\KM}{{\mathbb{K}}}

\newcommand{\cO}{{\mathcal O}}

\newcommand{\mm}{{\mathcal M}}
\newcommand{\FC}{{\mathcal F}}

\renewcommand{\a}{{\mathfrak a}}

\newcommand{\g}{{\mathfrak{g}}}
\newcommand{\h}{{\mathfrak{h}}}

\newcommand{\Uhg}{{U_\hbar({\g})}}

\newcommand{\cT}{{{}^cT}}
\newcommand{\uTTU}{{\underline{{}^cTTU'}}}
\newcommand{\CuTTU}{{C(\underline{{}^cTTU'})}}
\newcommand{\uTCuTTU}{{\underline{{}^cTC(\underline{{}^cTTU'})}}}

\newcommand{\uTA}{{\underline{{}^cTA}}}
\newcommand{\uTCuTA}{{\underline{{}^cTC(\underline{{}^cTA})}}}

\newcommand{\uTE}{{\underline{{}^cTE}}}
\newcommand{\uTpE}{{\underline{{}^cT(E)}}}
\newcommand{\uTpEu}{{\underline{{}^cT(E_1)}}}
\newcommand{\uTpEd}{{\underline{{}^cT(E_2)}}}

\newcommand{\uTCuTE}{{\underline{{}^cTC(\underline{{}^cTE})}}}

\newcommand{\dcH}{b_{\cH}}

\newcommand{\Dp}{D_{\poly}}
\newcommand{\Tpi}{T_{\poly}^{\inv}}

\newcommand{\Dpi}{D_{\poly}^{\inv}}
\newcommand{\fa}{\varphi_{\alg}}
\newcommand{\fg}{\varphi_{{\Ger}_\infty}}
\newcommand{\fl}{\varphi_{\Lie}}
\newcommand{\fld}{\varphi_{\Lie}'}
\newcommand{\fgi}{\varphi_{\Gi}}
\newcommand{\fgid}{\varphi_{\Gi}'}
\newcommand{\fgidu}{{{\varphi'}^1_{\Gi}}}

\begin{document}

\title[Quantization of coboundary $r$ matrices]
{Formality theorem for Lie bialgebras and quantization of
coboundary $r$-matrices}

\begin{abstract} 
Let $(\g,\delta_\hbar)$ be a Lie bialgebra.
Let $(\Uhg,\Delta_\hbar)$ a quantization of $(\g,\delta_\hbar)$ through 
Etingof-Kazhdan functor.
We prove the existence of a $L_\infty$-morphism
between the Lie algebra $C(\g)=\Lambda(g)$ and 
the tensor algebra $TU=T(\Uhg[-1])$ with
Lie algebra structure given by the Gerstenhaber bracket.
When $(\g,\delta_\hbar,r)$ is a coboundary Lie bialgebra,
we deduce from the formality morphism the
existence of a quantization $R$ of $r$.

\end{abstract}

\author{Gilles Halbout}
\address{IRMA (CNRS), rue Ren\'e Descartes, F-67084 Strasbourg, France}
\email{halbout@@math.u-strasbg.fr}

\maketitle

\centerline{0. {\textsc{Introduction}}}

\bigskip

Let $\KM$ be a field of characteristic $0$ and $\hbar$ a
formal parameter. Let $(\g,[-,-],\delta_\hbar)$ be a Lie bialgebra
over $\KM$ (all our objects will be $\KM[[\hbar]]$-modules). 
Using Etingof-Kazhdan quantization 
functor, one can construct a quantization 
$(\Uhg,\Delta_\hbar)$ of $(\g,\delta_\hbar)$. Let us denote $C(\g)=S(\g[-1])$
the free graded commutative algebra generated by $\g$;
$(C(\g),[-,-],\wedge,\delta_\hbar)$ is a differential Gerstenhaber algebra.
Let us also denote
 $TU=T(\Uhg[-1])$ the tensor algebra over $\Uhg$
(when $\delta_\hbar=0$, $U=U(\g)$, the enveloping algebra of $\g$).
More generally, we denote $TE=T(E[-1])$ the free tensor algebra of a
graded vector space $E$ and $\cT E=T(E[1])$ the cofree tensor coalgebra of $E$.
One can see the elements of $\g$ as invariant (under left action)
vector fields on the manifold $G$ where
$G$ is a connected group whose Lie algebra is $\g$.
In that framework, $C(\g)$ corresponds to the Gerstenhaber
algebra ${\Tpi}$ of invariant multivector fields on $G$
equipped with Schouten bracket.
The space $TU$ corresponds to the space ${\Dpi}$ of multidifferential
operators. The space $\Dp$ of Hochschild cochains carries a graded differential Lie algebra structure
when equipped with the Hochschild coboundary and
the Gerstenhaber bracket. Tamarkin proved in \cite{Tam98p} that the space $\Dp$
carries a $G_\infty$ structure.
In this paper, for general Lie bialgebra case, we prove:
\begin{theorem}
\label{theo:map}
\begin{enumerate}
\item There exists a $G_\infty$-structure on $TU$, whose underlying $L_\infty$-structure is the one given by the differential graded Lie structure
with deformed Gerstenhaber bracket and coHochschild differential.
\item There exists $\varphi$, a $L_\infty$-quasi-isomorphism between $C(\g)$ and $TU$ for 
the corresponding Lie algebra structures, such that the associated 
morphism of complex $\varphi^1$ maps $v \in C(\g)$ to its alternation $\Alt(v) \in TU \mod \hbar$.
\end{enumerate}
\end{theorem}
Definitions of $G_\infty$ and $L_\infty$-structures will be recalled
in section \ref{sec:structure} as well as the fact that $G_\infty$-algebras have canonical underlying $L_\infty$-structure.  This theorem generalises a result
of Calaque (\cite{C}) when $\delta_\hbar=0$ and answers to a conjecture of Tamarkin and
Tsygan (\cite{TT}).

\medskip

Suppose now that $(\g,r,Z)$ is a (finite-dimensional) coboundary Lie bialgebra over $\KM$. 
This means that $\g$ is a Lie bialgebra, 
the Lie cobracket $\delta_\hbar$ is the coboundary of an element
$r\in\Lambda^2(\g)$: 
$\delta_\hbar(x)=\hbar [x \otimes 1 + 1 \otimes x,r]$ for any $x \in \g$. 
This condition means that $Z:= \CYB(r)$ belongs to 
$\Lambda^3(\g)^\g$
(here $\CYB$ is the l.h.s. of the classical Yang-Baxter
equation). Quasi-triangular and triangular Lie bialgebras 
are particular cases of this definition. 

\medskip

Let us recall the definition of a quantization of $(\g,r,Z)$ in
coboundary Hopf algebra (\cite{Dr:QG}):
\begin{definition}
\label{defi:quanti}
A algebra $(U_\hbar(\g),R)$
is a coboundary quantization of $(\g,r,Z)$ if:
\begin{enumerate}
\item $(U_\hbar(\g),\Delta_\hbar)$ is a quantization 
of $\g$,
\item $R$ is an invertible element of
$U_\hbar(\g)^{\otimes 2}$ such that $R=1 + \hbar r + O(\hbar^2)$ that satisfies:
\item 
$R^{1,2}(\Delta_\hbar \otimes \id)(R) = R^{2,3}(\id \otimes \Delta_\hbar)(R)$,
\item $(\varepsilon \otimes \varepsilon)(R)=1,$
\item $R^{2,1}=R^{-1}.$
\item $R$ twists $\Delta_\hbar$ into $\Delta_\hbar^\op$:
$$R\Delta_\hbar(a)R^{-1}=\Delta_\hbar^{2,1}(a), \hskip1cm a \in U_\hbar(\g),$$
where $\Delta_\hbar^{2,1}=\Delta_\hbar^\op$ is the opposite comultiplication.
\end{enumerate}
\end{definition}
Here the notation $R^{i,j}$ corresponds to the coproduct-insertion map and will
be recalled at the end of this introduction. Using the formality map
of Theorem \ref{theo:map}, we prove the existence of a quantization of coboundary
$r$-matrices $r$.
\begin{theorem}
\label{theoprinc}
Let $(\g,r,Z)$ be a coboundary Lie bialgebra.
There exists a Hopf algebra $(U_\hbar(\g),$ $\Delta_\hbar)$ and an element
$R$ in
$U_\hbar(\g)^{\otimes 2}$ satisfying the first four properties
of Definition \ref{defi:quanti}.
\end{theorem}
\medskip

In section \ref{sec:structure}, we recall definitions of $G_\infty$
and $L_\infty$-structures. $G_\infty$-structures on a vector
space $A$ are defined on the cofree Lie coalgebra 
$S(\underline{\cT E}[1])$, where $\underline{\cT E}$ is the quotient of the cofree tensor coalgebra
$TE$ of a vector space $E$
by the image of the shuffles. 
We also recall the existence of two exact functors : 
${\LG}$, $\h \to C(\h)[1]$ between the categories of differential
Lie bialgebras and of differential Gerstenhaber algebras, viewed as
$G_\infty$-algebras and
$\LGi$, $\uTE \to E$ between the categories of differential
Lie bialgebras $\uTE$ and of $G_\infty$-algebras.

\medskip

In section \ref{sec:EK}, we recall Drinfeld duality between  QUE 
(Quantum Universal Enveloping) and QFSH 
(Quantum Formal Series Hopf) algebras.
In particular we recall the existence of functors
$(-)'$: QUE $\to$ QFSH and $(-)^\vee$: QFSH $\to$ QUE.
We then recall Etingof-Kazhdan quantization/dequantization functors.

\medskip

In section \ref{sec:brace}, we prove the existence of a bialgebra structure on $\cT TU$
coming from brace operations.
Moreover, we prove that
$\cT TU'$ is a QFSH algebra deforming the shuffle/cofree (for product/coproduct) structure on $\cT TU'$.
Then using Etingof-Kazhdan dequantization functor $\DQ$, we show that
the QUE algebra $(\cT TU')^\vee$  is sent
to $\uTTU$. So $TU'$ has a $G_\infty$-structure.
This structure reduces to a differential Lie algebra structure  with Gerstenhaber bracket and  coHochschild
differential defined using the coproduct $\Delta_\hbar$.
\medskip

In section \ref{sec:morphism}, we prove the existence of 
a bialgebra quasi-isomorphism $\fa$ : $U\to (\cT TU')^\vee$
and of a Lie bialgebra quasi-isomorphism $\fl$ : $\uTE \to \uTCuTE$ (for every vector space $E$).

\medskip

In section \ref{sec:inverse}, we recall the existence of an inverse map for any
$L_\infty$-quasi-isomorphism between $L_\infty$-algebras.

\medskip

In section \ref{sec:final}, we deduce the existence of a $L_\infty$-quasi-isomorphism between
$C(\g)$ and $TU$. The proof can be summarised in the following diagram:

\begin{equation*}
\begin{array}{ccccccccc}
\uTCuTTU & &\CuTTU[1] & = & \CuTTU[1]  &&\uTTU &&(\cT TU')^\vee\\[0.3cm]
\uparrow^{\fl} & \stackrel{{\LGi}}{\longrightarrow}& \uparrow^{\fgid} & & \uparrow^{\fg}&
\stackrel{{\LG}}{\longleftarrow}& \uparrow^{\fld}&\stackrel{{\DQ}}{\longleftarrow}& \uparrow^{\fa}\\[0.3cm]
\uTTU & &TU'[1]  &  & C(\g)[1]  &&\g &&U=\Uhg.
\end{array}
\end{equation*}
Thus the composition of $\fg$ with the inverse of $\fgi$ gives the wanted quasi-isomorphism.
From this, we prove Theorem \ref{theoprinc}.

\medskip

In the last section, we make some remarks on possible applications and
related open questions.
In particular, we prove that the coHochschild complex $(TU',\dcH)$ is quasi-isomorphic
to the complex $(C(\g),\delta_\hbar)$.

\medskip

\subsection*{Notations}

We use
the standard notation for the coproduct-insertion maps:
we say that an ordered set is a pair of a finite set $S$ and a 
bijection $\{1,\ldots,|S|\} \to S$. 
For $I_1,\dots,I_m$ disjoint ordered subsets of $\{1,\dots,n\}$,
$(U,\Delta)$ a Hopf algebra and $a \in  U^{\otimes m}$,
we define
$$a^{I_1,\dots,I_m}= \sigma_{I_1,\ldots,I_m} \circ 
(\Delta^{(|I_1|)}\otimes \cdots \otimes \Delta^{(|I_m|)})(a),
$$
with $\Delta^{(1)}=\id$, $\Delta^{(2)}=\Delta$, 
$\Delta^{(n+1)}=({\id}^{\otimes n-1} \otimes \Delta)\circ \Delta^{(n)}$, 
and $\sigma_{I_1,\ldots,I_m} : U^{\otimes \sum_i |I_i|} \to 
U^{\otimes n}$ is the 
morphism corresponding to the map $\{1,\ldots,\sum_i |I_i|\} 
\to \{1,\ldots,n\}$ taking $(1,\ldots,|I_1|)$ to $I_1$, 
$(|I_1| + 1,\ldots,|I_1| + |I_2|)$ to $I_2$, etc. 
When $U$ is cocommutative, this definition depends only on
the sets underlying $I_1,\ldots,I_m$.  

\smallskip

Until the end of this paper, although we will often omit to mention it, we will always deal with graded structures.

\medskip

\subsection*{Acknowledgements}
I would like to thank B. Enriquez and G. Ginot 
for discussions and D. Calaque and V. Dolgushev for
comments on this paper. I am very embedded to D. Tamarkin for his explanations.
I wish also to thank my cowriters of \cite{BGHHW} with whom I studied
 structures ``up to homotopy''.

\section{$G_\infty$-algebras}\label{sec:structure}

\medskip

\subsection{Definitions}

\medskip

Let us recall definitions of $L_\infty$-algebras and $L_\infty$-morphisms.
We have denoted $TA=T(A[-1])$ the free tensor algebra of $A$
which, equipped with the coshuffle coproduct, is a bialgebra.
We also have denoted $C(A)=S(A[-1])$ the free graded commutative algebra
generated by $A[-1]$, seen as a quotient
of $TA$. The coshuffle coproduct is still well defined on $C(A)$ which 
becomes a cofree cocommutative coalgebra on $A[-1]$.
We also denote $\Lambda A=S(A[1])$, the analogous graded commutative algebra
generated by $A[1]$ (in particular, for $A_1,A_2 \in A$, $A_1 \Lambda A_2$ stands for the
corresponding quotient of $A_1[1] \otimes A_2[1]$ in $\Lambda A$).
We will use the notations $T^nA$, $\Lambda^nA$ and $C^n(A)$ for the elements
of degree $n$.

\begin{definition}\label{DefinitionL_infini}
A vector space $A$ is endowed with a $L_\infty$-algebra 
(Lie algebra ``up to homotopy'') structure
if there are degree one linear  maps $d^{1,\dots,1}$: 
$\Lambda^kA \rightarrow A[1]$ such that 
the asociated coderivations 
(extended with respect to the
cofree cocommutative structure on $\Lambda A$) d: $\Lambda A \rightarrow \Lambda A$, 
satisfy $d \circ d=0$ where $d$ is the coderivation 
$$d=d^1+d^{1,1}+\cdots+d^{1,\dots,1}+\cdots.$$
\end{definition}
In particular, a differential Lie algebra $(A,b,[-,-])$ is
a $L_\infty$-algebra with structure maps $d^1=b[1]$, $d^{1,1}=[-,-][1]$ and
$d^{1,\dots,1}$: 
$\Lambda^k A \rightarrow A[1]$ are $0$ for $k \geq 3$.
One can now define the generalisation of Lie algebra morphisms:
\begin{definition}\label{Linfini}
A $L_{\infty}$-morphism between two $L_\infty$-algebras 
$(A_1,d_1=d^1_1+\cdots)$ and $(A_2,d_2=d^1_2+\cdots)$ is 
a morphism of codifferential cofree coalgebras, of degree $0$,
$$
\varphi: ~(\Lambda A_1,d_1) \rightarrow (\Lambda A_2,d_2).
$$
\end{definition}
\noindent In particular 
$\varphi \circ d_1 = d_2 \circ \varphi $. 
As $\varphi$ is  a morphism of cofree cocommutative  coalgebras, $\varphi$
is determined by its image on the cogenerators, i.e., by its components:
$\varphi^{1,\dots,1}: \Lambda^k A_1\to A_2[1].$ 

\medskip

Let us denote $\cT( E)$ the cofree tensor coalgebra of $E$
(with coproduct $\Delta'$). Equipped with the
shuffle product $\bullet$ (defined on the cogenerators  $\cT (E) \otimes  \cT (E) \to E$
as $\pr \otimes \varepsilon + \varepsilon \otimes \pr$, where $\pr$ :  $\cT (E) \to E$ is
the projection and $\varepsilon$ is the counit), it is a bialgebra. 
Let $\cT (E)^+$ be the augmentation ideal.
We denote $\uTpE= \cT (E)^+/(\cT (E)^+ \bullet \cT (E)^+)$ the quotient by the shuffles.
It has a graded cofree Lie coalgebra structure (with coproduct $\delta = \Delta' -{\Delta'}^\op$).
Then $S(\uTpE[1])$ has a structure of cofree coGerstenhaber algebra (i.e. equipped with
cofree coLie and cofree cocommutative coproducts satisfying compatibility condition).
We use the notation $\underline{\cT^n (E)}$ for the elements of degree $n$.

\begin{definition}\label{DefinitionG_infini}
A vector space $E$ is endowed with a $G_\infty$-algebra 
(Gerstenhaber algebra ``up to homotopy'') structure
if there are degree one linear  maps $d^{p_1,\dots,p_k}$: 
$\underline{\cT^{p_1}(E)} \Lambda \cdots \Lambda \underline{\cT^{p_k}(E)}
\subset \Lambda^k \uTE \to E[1]$
such that the associated coderivations 
(extended with respect to the
cofree coGerstenhaber structure on $\Lambda \uTpE$) d:
$\Lambda \uTpE \rightarrow \Lambda \uTpE$
satisfies  $d \circ d=0$ where $d$ is the coderivation 
$$d=d^1+d^{1,1}+\cdots+d^{p_1,\dots,p_k}+\cdots.$$
\end{definition}
In particular we have
\begin{remark}
\label{rem:shift}
If $(E,b,[-,-],\wedge)$ is a differential Gerstenhaber algebra, then $E[1]$  is
a $G_\infty$-algebra with structure maps $d^1=b[1]$, $d^{1,1}=[-,-][1]$, $d^2=\wedge[1]$ and
other $d^{p_1,\dots,p_k}$: 
${{\underline{{}^cT^{p_1}(E[1])}}}$ $ \Lambda \cdots \Lambda {{\underline{{}^cT^{p_k}(E[1])}}} \rightarrow E[2]$ are $0$.
\end{remark}
One can finally define the generalisation of Gerstenhaber algebra morphisms:
\begin{definition}\label{Ginfini}
A $G_{\infty}$-morphism between two $G_\infty$-algebras 
$(E_1,d_1=d^1_1+d_1^2+\cdots)$ and $(E_2,d_2=d^1_2+d_2^2+\cdots)$ is 
a morphism of differential coGerstenhaber coalgebras, of degree $0$,
$$
\varphi: ~(\Lambda \uTpEu,d_1) \rightarrow (\Lambda\uTpEd,d_2).
$$
\end{definition}
\noindent 
In particular $\varphi \circ d_1 = d_2 \circ \varphi $. 
As $\varphi$ is  a morphism of cofree coGerstenhaber coalgebras, $\varphi$
is determined by its image on the cogenerators, i.e., by its components:
$\varphi^{p_1,\dots,p_k}$: 
$\underline{\cT^{p_1}(E_1)} \Lambda 
\cdots $ $\Lambda \underline{\cT^{p_k}(E_1)}\to E_2[1].$ 

\smallskip

\subsection{Functors ${\LG}$ and $\LGi$}

\medskip

Let $(\h,\delta,d)$ be a differential Lie bialgebra. Let $C(\h)=S(\h[-1])$ be the free graded
commutative algebra generated by $\h$. Recall from the previous subsection that $C(\h)$ is also a cofree coalgebra and that 
coderivations $C(\h) \to C(\h)$ are defined by their images in $\h$. Thus, one
easily checks that the coderivation $[-,-]$: $C(\h) \to C(\h)$ extending the Lie
bracket (with degree shifted by one) defines a Lie algebra structure
on $C(\h)$ and that $(C(\h),[-,-],\wedge)$ is a Gerstenhaber algebra,
where $\wedge$ is the commutative product: $(\alpha,\beta) \mapsto
\alpha \wedge \beta= \alpha \Lambda \beta$
on $C(\h)$. Moreover,
one can extend maps $d$: $\h \to \h$ and $\delta$: $\h \to S^2(\h[-1])$ on
the free commutative algebra $C(\h)$ so that
$(C(\h),[-,-],\wedge,d+\delta)$ is a differential Gerstenhaber algebra. In fact, the
differential $\delta$ is the Chevalley Eilenberg differential: the space
$C(\h)=S^*(\h[-1])$ is isomorphic to the standard complex $(\Lambda^*(\h))[-*]$ and
$\delta$ is simply the differential given by 
the underlying Lie algebra structure of $\h$.
Moreover any morphism $f:\h_1\to \h_2$ can be 
extended into a morphism $C(f):C(\h_1) \to C(\h_2)$ 
of free commutative algebras thanks to the inclusion 
$h_2\subset C(\h_2)$. This morphism is easily seen to be  
a differential Gerstenhaber algebra morphism.
Thus, we have defined a functor ${\LGG}$ 
from differential Lie bialgebras to differential Gerstenhaber algebras which sends 
$\h$ to $C(\h)$.
 This functor  is exact.
 As the differential $\delta$ and $d$ anticommutes by construction, 
 the complex $(C(\h),d,\delta)$ is a (first quadrant) bicomplex.
  Hence a quasi-isomorphism $(\h_1,d_1)\to (\h_2,d_2)$ 
  induces a quasi-isomorphism 
  $(C(\h_1),d_1,\delta_1)\to (C(\h_2),d_2,\delta_2)$.
 Then, thanks to Remark \ref{rem:shift}, we get a functor $\LG$: 
 $\h \to C(\h)[1]$ from 
  differential Lie bialgebras to $G_\infty$-algebras.

\medskip

Let us now define the functor $\LGi$. Consider the category ${\rm CFDLB}$ of  
differential Lie bialgebras which are cofree as a Lie coalgebra. 
In other words we are interested in cofree Lie coalgebra $\uTpE$ 
on a graded vector space $E$ together with a differential $\ell$ 
and a cobracket $L$ on $\uTpE $ that makes it a differential Lie bialgebra. 
As $\uTpE$ is cofree, the differential is uniquely determined 
by its restriction to cogenerators $l^p$: $\underline{\cT^p(E)} \to E$. 
Similarly, the Lie bracket is uniquely determined by maps $L^{p_1,p_2}$:  
$\underline{\cT^{p_1}(E)}\Lambda \underline{\cT^{p_2}(E)} \to E$.
Now, this data determines on $E$  a structure of $G_\infty$-algebra
with structure maps given 
by $d^{p_1,\dots,p_k}$: 
$\underline{\cT^{p_1}(E)} \Lambda \cdots \Lambda \underline{\cT^{p_k}(E)}
\to E$, with $d^{p_1,\dots,p_k}=0$ for $k >2$ and $d^{p_1,p_2}=L^{p_1,p_2}$
and $d^{p}=l^{p}$ (with degrees shifted by one). In fact, 
according to Definition \ref{DefinitionG_infini}, 
a $G_\infty$-structure on $E$ is given by a differential 
$d$ on $\Lambda \uTpE$ which as a space is isomorphic 
to the standard Chevalley-Eilenberg complex of the differential Lie algebra 
$(\uTpE,\ell,L)$. The differential defined above is simply 
the Chevalley-Eilenberg differential. In particular $d$ is the sum $d=d^1+d^2$ 
with $d^1=\sum_{p\geq 1} l^p$ and $ d^2=\sum_{p_1,p_2\geq 1}L^{p_1,p_2} $ 
and $(\Lambda \uTpE, d^1,d^2)$ is a bicomplex. 
Moreover, a morphism $\varphi$: $\underline{T(E_1)} \to \underline{T(E_2)}$ 
of differential Lie bialgebras is determined by its 
restriction to cogenerator of the cofree Lie coalgebra structure, 
that is to say by maps: $\varphi^n$: $\underline{T^n(E_1)} \to E_2$. 
It determines a $G_\infty$-morphism
$E_1 \to E_2$ (with the $G_\infty$-structures defined above) 
defined by maps $\varphi^n$: $\underline{T^n(E_1)} \to E_2$,
other being $0$. This is simply the functoriality of 
the Chevalley-Eilenberg complex. Thus we have defined a 
  functor from ${\rm CFDLB}$ to the category of $G_\infty$-algebras. 
 This functor is exact. A quasi-isomorphism of differential Lie bialgebras 
 $(\uTpEu,\ell_1)\to  (\uTpEd,\ell_2)$ induces a quasi-isomorphism 
 $(\Lambda \uTpEu, d^1_1)\to (\Lambda \uTpEd, d^1_2)$, hence, 
 as  $(\Lambda \uTpE, d^1,d^2)$ is a (first quadrant) bicomplex,  
 a quasi-isomorphism $(\Lambda \uTpEu, d_1)\to (\Lambda \uTpEd, d_2)$.

\medskip

Untill the end of the paper, we will use the notations $TE$ for $T(E[-1])$ 
and $\uTE$ 
for $\underline{T(E[1])}$.

\medskip

\section{Etingof-Kazhdan functors}\label{sec:EK}

\medskip

\subsection{Duality of QUE and QFSH algebras} \label{subsec:duality}

\medskip

we recall some facts from \cite{Dr:QG} (proofs can be found in
\cite{Gav}). Let us denote by ${\bf QUE}$ the category of quantized universal
enveloping (QUE) algebras and by ${\bf QFSH}$ the  category of quantized formal
series Hopf (QFSH) algebras. We denote by ${\bf QUE}_{\fd}$ and
${\bf QFSH}_{\fd}$ the subcategories corresponding to finite dimensional 
Lie bialgebras. 

We have contravariant functors ${\bf QUE}_{\fd} \to {\bf QFSH}_{\fd}$, 
$U\mapsto U^*$ and ${\bf QFSH}_{\fd} \to {\bf QUE}_{\fd}$, 
$\cO\mapsto \cO^\circ$. These functors are inverse to each other. 
$U^*$ is the full topological dual of $U$, i.e., the space of all  
continuous (for the $\hbar$-adic topology) $\KM[[\hbar]]$-linear maps
$U \to \KM[[\hbar]]$. 
$\cO^\circ$ is the space of continuous $\KM[[\hbar]]$-linear forms 
$\cO\to \KM[[\hbar]]$, 
where $\cO$ is equipped with the $\mm$-adic topology (here $\mm\subset \cO$
is the maximal ideal). 

We also have covariant functors ${\bf QUE} \to {\bf QFSH}$, $U\mapsto U'$
and ${\bf QFSH} \to {\bf QUE}$, $\cO\mapsto \cO^\vee$. These functors are
also inverse to each other. 

$U'$ is a subalgebra of $U$
defined as follows:
for any ordered subset 
$ \, \Sigma = \{i_1, \dots, i_k\} \subseteq \{1, \dots, n\} \, $ 
with  $ \, i_1 < \dots < i_k \, $,  \, define the morphism 
$ \; j_{\scriptscriptstyle \Sigma} : U^{\otimes k} \longrightarrow
U^{\otimes n} \; $  by  $ \; j_{\scriptscriptstyle \Sigma}
(a_1 \otimes \cdots \otimes a_k) := b_1 \otimes \cdots \otimes
b_n \; $  with  $ \, b_i := 1 \, $  if  $ \, i \notin \Sigma \, $  and 
$ \, b_{i_m} := a_m \, $  for  $ \, 1 \leq m \leq k $;  then set 
$ \; \Delta_\Sigma := j_{\scriptscriptstyle \Sigma} \circ \Delta^{(k)}
\, $,  $ \, \Delta_\emptyset := \Delta^{(0)} \, $,  and
 $$  \delta_\Sigma := \sum_{\Sigma' \subset \Sigma} {(-1)}^{n- \left|
\Sigma' \right|} \Delta_{\Sigma'}  \; ,  \qquad \quad  \delta_\emptyset
:= \epsilon \; .   $$  
We shall also use the notation  $ \, \delta^{(n)} := \delta_{\{1, 2,
\dots, n\}} \, $,  $ \, \delta^{(0)}
:= \delta_\emptyset \, $,  and the useful formula  
  $$  \delta^{(n)} = {({id}_{\scriptscriptstyle U} - \epsilon)}^{\otimes n}
\circ \Delta^{(n)} \, .  $$  
Finally, we define
  $$  U' := \big\{\, a \in U \,\big\vert\, \delta^{(n)}(a) \in
h^n U^{\otimes n} \, \big\} \quad ( \subseteq U )  $$
and endow it with the induced topology.  

\medskip

On the other way, $\cO^\vee$
is the $\hbar$-adic completion of $\sum_{k\geq 0} \hbar^{-k} \mm^k \subset 
\cO[1/\hbar]$.

\medskip

We also have canonical isomorphisms $(U')^\circ \simeq (U^*)^\vee$
and $(\cO^\vee)^* \simeq (\cO^\circ)'$. 

\medskip

\subsection{The functor $\DQ$}

\medskip

In \cite{GH}, a generalisation of Etingof-Kazhdan theorem (\cite{EK}) was proven
in an appendix by Enriquez and Etingof:
\begin{theorem}
\label{thm:dqdg}
We have an equivalence of categories
$$DQ_\Phi~:~\DGQUE \to \DGLBA_h$$
from the category of differential graded  
quantized universal enveloping super-algebras to that of differential graded Lie
super-bialgebras such that if
$U \in \Ob(\DGQUE)$ and $\a=DQ(U)$, then
$U/hU=U(\a/h\a)$, where $U$ is the universal algebra functor, taking a
differential graded Lie super-algebra to a differential graded super-Hopf 
algebra.
\end{theorem}
Here $\Phi$ is a Drinfeld assoicator. We will use any of these functor and denote it $\DQ$.

\medskip

\section{Bialgebra structure on $\cT TU$}\label{sec:brace}

\medskip

\subsection{Brace structures}
Here, we will define a structure on $\cT TU$.
Following \cite{br1,br2,br3}, it is a brace structure (although
usually brace structures are defined on the space of Hochschild co\-chains of 
an algebra) but we will not recall the definition of brace structure as we will
only need the result of the following subsection. 
Let us define those structures for a general Hopf algebra.
More precisely, let $(U,\Delta_\hbar,\times)$ be a Hopf algebra (in our case $U$ will be
the Etingof-Kazhdan quantization $\Uhg$ of the Lie bialgebra $\g$).
We will define a brace structure on the cofree tensor coalgebra $\cT TU$ of
the free tensor algebra $T(U[-1])$. 
To distinguish the two tensor products, we denote $\otimes$ the
tensor product on $TU$ and $\boxtimes$ the tensor product
on $\cT TU$.
\begin{definition}
\label{brace}
We define a brace structure on $\cT TU$ by extending the following maps
given on the cogenerators of the cofree coalgebra $\cT TU$:
\begin{enumerate}
\item $B^0=0$,
\item $B^1=\dcH$ (the coHochschild coboundary on $TU$),
\item $B^2$ : $\alpha \boxtimes \beta \mapsto \alpha \otimes \beta$,
\item $B^n=0$  for $n >2$,
\item $B^{0,1}=B^{1,0}=\id$,
\item $B^{0,n}=B^{n,0}=0$ for $n >1$,
\item $B^{1,n}$ : $(\alpha,\beta_1\boxtimes \cdots \boxtimes \beta_n) \mapsto$ 
\begin{multline*}
\sum_{\stackrel{0 \leq i_1,\dots,i_m+k_m\leq n}
{i_l+k_l\leq i_{l+1}}} (-1)^\varepsilon \alpha^{1,\dots,i_1+1\cdots i_1+k_1,\dots,i_m+1 \cdots i_m+k_m,
\dots,n} \times \\
1^{\otimes i_1} \otimes \beta_1 \otimes 1^{\otimes i_2-(i_1+k_1)}\otimes \beta_2 \otimes \dots \otimes \beta_n \otimes
1^{\otimes n-(i_m+k_m)},
\end{multline*}
where $k_s=|\beta_s|$, $n=|\alpha|+\sum_sk_s - m$ and $\varepsilon = \sum_s (k_s-1)i_s$,
\item $B^{k,l}=0$ for $k >1$.
\end{enumerate}
\end{definition}
Note that, when $U=U(\g)$, the enveloping
algebra of a Lie algebra $\g$, $T(U[-1])$ can been seen as the space of
invariant multidifferential operators over the Lie group corresponding to $\g$
and in that case, our definition coincides with those of \cite{br1,br2,br3}.

\medskip

\subsection{From Hopf algebra $U$ to Bialgebra structure on $\cT TU$}

\medskip

Still following \cite{br1,br2,br3}, we have:
\begin{theorem}
The brace structure of Definition \ref{brace} defines a differential bialgebra
structure on the cofree tensor coalgebra $\cT TU$, with product  $\star$
extending $\sum B^{p_1,p_2}$ and differential $d$ extending
$\sum B^p$.
\end{theorem}
{\em Proof.}
To prove the associativity of $\star$, one has to check the following equation
for $\alpha,\beta_1,\dots, \beta_l,$ $\gamma_1,\dots,\gamma_l \in TU$:
\begin{multline*}
\sum_{0 \leq i_1 \leq \cdots \leq i_l \leq m}(-1)^\varepsilon
\{\alpha\}\{\gamma_1,\dots,\gamma_{i_1},\{\beta_1\}  \{ \gamma_{i_1+1},\dots \}  ,
\dots, 
 \gamma_{i_l} ,\{ \beta_{l}\}\{  \gamma_{i_l+1} ,\dots\},\dots,\gamma_{m}\} \cr
 =\{\{\alpha\}\{\beta_1,\dots,\beta_l\}\}\{\gamma_1,\dots,\gamma_m\},
\end{multline*}
where $\varepsilon=\sum_{p=1}^l(|\beta_p|-1)\sum_{q=1}^{i_p}(|\gamma_q|-1)$
and $\{\alpha\}\{\beta_1,\dots,\beta_l\}$ is
$B^{1,l}(\alpha,\beta_1 \boxtimes \cdots \boxtimes \beta_l)$.
This equation is a consequence of the associativity, coassociativity 
and compatibility of
$\times$ and $\Delta$:
\begin{align*}
(\Delta_\hbar^{(k)}\otimes \id^{\otimes l} \otimes \Delta_\hbar^{(m)})(\alpha)&= 
(\Delta_\hbar^{(k)} \otimes \id^{\otimes l+m}) \circ (\id^{\otimes l+1}
 \otimes \Delta_\hbar^{(m)})(\alpha)\cr
\hbox{ and }& \cr
(\Delta_\hbar^{(k)}\otimes   \id^{\otimes n+p-2}  )((\Delta_\hbar^{(n)}\otimes  
\id^{\otimes p-1}   )(\alpha) \times (\beta \otimes 1^{\otimes p-1}       ))&= 
(\Delta_\hbar^{(n+k-1)} \otimes  \id^{\otimes p-1}   )(\alpha) \cr
&\times (\Delta_\hbar^{(k)} \otimes  \id^{\otimes n+p-2}     )
(\beta \otimes  1^{\otimes p-1}),
\end{align*}
For $\alpha,\beta\in  TU$ of degree $p$ and $n$.

\smallskip

One can then notice that  the map $d$ is the
commutator, with respect to the product $\star$, of the element $1 \otimes 1 \in TU \subset \cT TU$.
Thus $d$ is compatible with the multiplication.
Finally, $d$ is a differential as $[1\otimes 1, 1 \otimes 1]_\star =0$.
\hfill \qed \medskip
\begin{remark}
Note that we have only used here the bialgebra structure of $U$. So one would
get a similar result when replacing $U$ with $U'$ (see section \ref{sec:EK}).
\end{remark}

\medskip

We have now:
\begin{proposition}
The algebra $\cT TU'$ is a QFSHA.
\end{proposition}
{\em Proof.}
This is known when $U=U(\hbar \g)$ (\cite{Tam98p,TT}  or \cite{GH}). Proof can be done in the
same way for any QUE algebra $\Uhg$: one considers the dual bialgebra $(\cT TU')^*$.
It is a free algebra and so a QUE algebra: as a $\KM[[\hbar]]$-module, it
is isomorphic to the enveloping algebra of the corresponding free Lie algebra
$(\uTTU')^*$. So $\cT TU'$ is a QFSH algebra
\hfill \qed \medskip

\begin{remark}
In this proof, we have shown that the corresponding differential Lie bialgebra to $(\cT TU)^\vee$ through
Etingof-Kazhdan dequantization functor $\DQ$ is isomorphic to
$\uTTU$ as a  $\KM[[\hbar]]$-module.
\end{remark}
Thus $TU'$ is  equipped with a $G_\infty$-structure (see section \ref{sec:structure}).

\medskip

Finaly, one easily checks that 
\begin{proposition}
The underlying differential Lie bialgebra structure
on $TU'$ corresponding to the differential Lie bialgebra structure on $\uTTU$
 is given by Gerstenhaber bracket 
$$[\alpha,\beta]_G=B^{1,1}(\alpha,\beta) - (-1)^{(|\alpha|-1)(|\beta|-1)}B^{1,1}(\beta,\alpha) $$
and coHochschild differential
$$\dcH(\alpha)=[1\otimes 1, \alpha]_G.$$
\end{proposition}

\medskip

\section{Bialgebra quasi-isomorphisms}\label{sec:morphism}

\medskip

\subsection{A bialgebra quasi-isomorphism $\fa$ : $U\to (\cT TU)^\vee$}

\medskip

Let us first define a bialgebra quasi-isomorphism
$\fa$ : $U'\to \cT TU'$ between the
bialgebra $(U',\Delta_\hbar,\times)$ (in our case $U=\Uhg$ the
Etingof-Kazhdan quantization of $\g$) and
the bialgebra $(\cT TU',\Delta,\star)$ whose
structure was described in the previous
section. The definition
of $U'$ was given in section \ref{sec:EK}. As $\cT TU'$ is a cofree coalgebra, as a coalgebra map $\fa$ is uniquely determined by its restriction $U\to TU'$ to cogenerators of $\cT TU'$. We define
$$\fa= \inc - \varepsilon 1,$$
where $\inc$ is actually the inclusion $U' \subset TU'$.
Equivalently, following \cite{Ta2},
let us consider $\mm$ the augmentation ideal of $U'$ and denote
\begin{equation*}
\begin{array}{rl}
\delta^{(2)}: \mm &\to \mm \otimes \mm,\\
x &\mapsto \Delta_\hbar(x)-(1 \otimes x + x \otimes 1).
\end{array}
\end{equation*}
We still use the notation $T{\mm}$ for $T(\mm[-1])$ and $\cT T\mm$ for
$\cT (T\mm[1])$. One now defines 
\begin{equation*}
\begin{array}{rl}
\bar{\delta}: \mm &\to T(\mm),\\
x &\mapsto x + \sum_{k \ge   2}\delta^{(k)}(x),
\end{array}
\end{equation*}
where $\delta^{(k)}$ is the $k-1$-th iterate of $\delta^{(2)}$ (we will set $\delta^{(1)}=\inc$).
This map is well defined thanks to the definition of $U'$.
Now inclusion $i$: $U'[-1] \to T(U'[-1])$ defines, after prolongation on the
cofree tensor coalgebra, a map $Ti$:
$T(U'[-1][1]) \to \cT TU'$. We have $\fa= Ti \circ \delta$.
We can then take into account the unit and define the map $\fa$ on 
$U'$: $\KM \oplus \mm \to \KM \oplus \cT T\mm$. 

As every maps are coalgebra morphisms, 
to check that $\fa$ is an algebra
morphism, one only has to check that the following diagram commutes after 
canonical projection 
$Pr$: $\cT TU' \to TU'[1]$:
\begin{equation*}
\begin{array}{ccccc}
U' \otimes U' & \stackrel{{\fa \otimes \fa}}{\longrightarrow} &\cT TU' \otimes \cT TU' &  & \\[0.3cm]
\downarrow_{\times} && \downarrow{\star} & & \\[0.3cm]
U'& \stackrel{{\fa }}{\longrightarrow}  &\cT TU' & \stackrel{{Pr}}{\longrightarrow}
&TU'[1].
\end{array}
\end{equation*}
Let us check the commutation for $\mu\otimes \eta \in \mm \otimes \mm$.
Commutation is obvious if one of the elements $\mu$ or $\eta$ is in $\KM$, so let us
consider the image of $\mu\otimes \eta \in \mm \otimes \mm$ through the
two paths of this diagram. 
Using the notations of section \ref{sec:brace}, we have
$\delta\alpha =\sum_n \delta^{(n)} \alpha = \sum \alpha_{n_1} 
\boxtimes \cdots \boxtimes \alpha_{n_k}$ and
$\delta\beta = \sum \beta_{m_1} 
\boxtimes \cdots \boxtimes \beta_{m_l}$
with $\alpha_{n_i},\beta_{m_j} \in U$. Note
that $B^{1,n}(\alpha_{n_i},\beta_{m_1} 
\boxtimes \cdots \boxtimes \beta_{m_l})$ is
$0$ for $n >1$ ($\alpha_{n_i}$ and the $\beta_j$'s are in  $\in U'$)
and is $\alpha_{n_i} \times \beta_{m_1}$ otherwise
and that $B^{p>1,q}=0$. So we get
$$Pr(\fa(\alpha) \star \fa(\beta))=\Pr (\alpha \times \fa(\beta))=\alpha \times \beta,$$
which is the commutation property.

\smallskip

Let us show now that $\fa$ is a quasi-isomorphism of complexes.
Recall that the differential $d$ on $\cT TU'$ is defined as extension of 
$B^1 + B^2$ (c.f. definition \ref{brace}).
Let us prove the following lemma
\begin{lemma}
We have $H_\cdot(\cT T(U'),d) \simeq U'$ in degree $1$ and is $0$ for degree $\geq 2$.
\end{lemma}
{\em Proof.}
\smallskip
The differential $d$ is the sum $d=d^1+d^2$ where $d^1$ and $d^2$ correspond respectively
to the maps $B^1$ and $B^2$. Thus $\cT TU'$ has a structure of bicomplexe.
Let us compute the homology with respect to $d^2$: 
$(\cT TU' \simeq \oplus \cT^n TU',d^2)$ is the Hochschild complexe associated
to $TU'$ and $TU'$ is a free associative algebra then
$$H_\cdot(\cT T(U'),d^2)\simeq H_\cdot(\cT^n TU',d^2)\simeq U',$$
concentrated in degree $1$. So $H_\cdot(\cT T(U'),d^1+d^2) \simeq \Ker(U',d^1)=U'$
as $d^1(H_\cdot(TU'))\cong 0$.
\hfill \qed \medskip 

\smallskip

Finally, we check that $\fa$ is a morphism of complexes. As before, it is enough to
check on the cogenerators that $Pr( d(\fa(\alpha)))=0$ for $\alpha \in U'$.
Still writing $\delta\alpha = \sum \alpha_{n_1} 
\boxtimes \cdots \boxtimes \alpha_{n_k}$ we get
$$Pr(d(\fa(\alpha)))=Pr(d^1(\alpha)+d^2(Ti\delta^{(2)}(\alpha)))=\dcH (\alpha)-\delta^{(2)}(\alpha)=0,$$
where $\dcH$ is the coHochschild codifferential.

\medskip

Thus we have a bialgebra quasi-isomorphism
$\fa$ : $U'\to \cT TU'$ . Applying to it the Drinfeld functor $(-)^\vee$, we get 
a bialgebra quasi-isomorphism
$\fa$ : $U\to (\cT TU')^\vee$.

\medskip

\subsection{A Lie bialgebra quasi-isomorphism $\fl$ : $\uTA \to \uTCuTA$}

\medskip

Let $A$ be a vector space.
Suppose now that the cofree Lie coalgebra $\uTA$
has a structure $(\uTA,\delta,[-,-],d)$ of a differential Lie bialgebra.
Using the functor ${\LG}$ (see section \ref{sec:structure}), one gets a differential Gerstenhaber algebra
$(C(\uTA),[-,-],\wedge,d+\delta)$.
One can extend the structure maps on
the cofree Lie coalgebra $\uTCuTA$
and one gets 
a differential cofree Lie bialgebra 
$(\uTCuTA,\delta',[-,-],d+\delta+\wedge)$ (we will set $d^1=d + \delta$ and  $d^2=\wedge$).
We can now prove the existence of a differential Lie bialgebra quasi-isomorphism $\fl$ : $\uTA \to \uTCuTA$.
Let us extand the inclusion $\inc$: $\uTA \to C(\uTA)$ to a coderivation
$\fl$: $\uTA \to \uTCuTA$ on the cofree Lie coalgebra. More explicitly, 
using the Lie cobracket $\delta$ of  $\uTA$, on defines
\begin{equation*}
\begin{array}{rl}
\bar{\delta}: \uTA &\to S(\uTA),\\
x & \mapsto \sum_{k \ge   1}\delta^{(k)}(x),
\end{array}
\end{equation*}
where $\delta^{(k)}$ is the $k-1$-th iterate of $\delta$ and $\delta^{(1)}=\inc $.
Now inclusion $i$: $\uTA[-1] \to C(\uTA)$ defines, after prolongation on the
cofree tensor Lie coalgebra, a map $Ti$:
$S(\uTA) \to \uTCuTA$. We have $\fl= Ti \circ \delta$.
We can reproduce the proof of the previous subsection (all structures
are much more simpler) and easily check that
$\fl$ : $\uTA \to \uTCuTA$ is a differential Lie bialgebra quasi-isomorphism.
 It is obviously a Lie bialgebra morphism. Moreover, for $\alpha \in \uTA$,
we get after projection on the cogenerators,
$$(d^1+d^2)(\fa(\alpha))=d^1(\alpha) + d^2(Ti \delta(\alpha))
=d(\alpha) -\delta(\alpha) + \delta(\alpha) = d(\alpha).$$
The fact that $\fl$ is a quasi-isomorphism can be proved as in the previous section.

\medskip

\section{Inversion of formality morphisms}\label{sec:inverse}

\medskip

Let us recall Theorem 4.4 of Kontsevich (\cite{K}):
\begin{theorem}
Let $\g_1$ and $\g_2$ be two $L_\infty$-algebras and $\FC$
be a $L_\infty$-morphism from $\g_1$ to $\g_2$.
Assume that
$\FC$ is a quasi-isomorphism. Then there exists an $L_\infty$-morphism
from $\g_2$ to $\g_1$ inducing the inverse isomorphism
between associated cohomology of complexes.
\end{theorem}

\begin{remark}
\label{greg}
We know from private communications the existence of a similar
$G_\infty$-version of this theorem. This result would imply the
existence of a corresponding $G_\infty$-morphism in Theorem \ref{theo:map}.
\end{remark}

\medskip

\section{Proof of the main theorems}\label{sec:final}

\medskip

\subsection{Proof of Theorem \ref{theo:map}}
 
\medskip

Let $(\g,\delta_\hbar)$ be a Lie bialgebra. We write
$\delta_\hbar=\hbar \delta_1 + \hbar^2 \delta_2 + \cdots$.
Let $(\Uhg,\Delta_\hbar)$ be the Etingof-Kazhdan canonical quantization
of $(\g,\delta_\hbar)$. 
We denote $U=\Uhg$ for short.
In section \ref{sec:brace}, we proved the existence of a bialgebra structure on $\cT TU$ and
thanks to section \ref{sec:morphism}, we have a
bialgebra quasi-isomorphism $\fa$: $U\to (\cT TU')^\vee$.
Thanks to Etingof-Kazhdan dequantization functor (see section
\ref{sec:EK}), and the fact
that $(\cT TU')^\vee$ is a QUE algebra quantizing $\uTTU$
(see section \ref{sec:brace}),
we get a Lie bialgebra quasi-isomorphism $\fld$: $\g\to\uTTU$,
a differential Lie bialgebra.
Using the exact functor ${\LG}$ (see section \ref{sec:structure}), we get a 
quasi-isomorphism of differential Gerstenhaber algebras $\fg$:
$C(\g)[1] \to \CuTTU[1] $.

\smallskip

According to section \ref{sec:morphism}, we also have a
differential Lie bialgebra quasi-isomorphism $\fl$ : $\uTTU \to \uTCuTTU$.
This quasi-isomorphism is sent to a $G_\infty$-quasi-isomorphism
$\fgid$: $TU'[1] \to \CuTTU[1]$ using the functor $\LGi$ defined in section \ref{sec:structure}.

\smallskip

Finally, in section \ref{sec:inverse}, we recalled the existence of an inverse map for any
$L_\infty$-quasi-isomorphism between $L_\infty$-algebras so
the $L_\infty$-restriction of the 
quasi-isomorphism $\fgid$ is invertible. Then one can define
$\varphi$: $C(\g) \to TU'$ as the composition of the $L_\infty$-restriction of $\fg$
with the inverse of the $L_\infty$-restriction of $\fgid$.

\medskip

One has to show now that $\varphi^1$ maps $v \in C(\g)$ to $\Alt(v) \in TU' \mod \hbar$.
Let us replace $\g$ with $\hbar \g$ in the previous construction.
Let $\hbar v $ be an element of $\hbar \g$.
Let us still call $\hbar v$ its image in the quantization $U$.
By definition of $\fa$, $\fa(\hbar v)=\hbar v \mod \hbar^2$ (once again,
we still use the same notation for its image in $(\cT TU')^\vee$.
So $\fld(\hbar v)=\hbar v \mod \hbar^2$ ($\hbar v$ on the right hand side is the
image of $\hbar v$ in $\uTTU$). Note that here appeared highly non trivial
terms in $O(\hbar^2)$. So $\fg(\hbar v_1 \wedge \cdots  \wedge \hbar v_n) = \hbar^n \Alt(v_1\otimes
\cdots \otimes v_n) \mod \hbar^{n+1}$.

Moreover, it is clear by definition, that $\fl$: $\uTTU \to \uTCuTTU$ is 
the identity map $\id \mod \hbar$
when restricted to $TU'$ thus so is 
the corresponding map of complexes ${\fgidu}$: $TU' \to \CuTTU$. Then $\varphi^1$:
$C^n(\hbar \g) \to TU'$ is $\Alt \mod \hbar^{n+1}$.
\hfill \qed \medskip 

\medskip

\subsection{Proof of Theorem \ref{theoprinc}}

\medskip

Suppose now that $(\g,\delta_\hbar,[-,-],r,Z)$ is a (finite-dimensional) coboundary Lie bialgebra over $\KM$. 
This means that
the Lie cobracket $\delta_\hbar$  is the coboundary of 
$r\in\Lambda^2(\g)$: 
$\delta_\hbar(x)=\hbar [x \otimes 1 + 1 \otimes x,r]$ for any $x \in \g$. 
Let $r'=-\hbar r$. This means that $ r'$ is a Mauer-Cartan element in $\hbar \g$:
$[r',r'] + \delta_\hbar(r')=0$. Let us write the $L_\infty$-morphism $\varphi$
of Theorem \ref{theo:map}: $\varphi= \sum_{k\geq 1} \varphi^{1,\dots,1}$,
where $ \varphi^{1,\dots,1}$:
$\Lambda^n \g \to TU'$.
Let us define $R'= 1 \otimes 1 + \sum_{k \geq 1} \frac{1}{n!}(\Lambda^n r')$.
By definition of $L_\infty$-morphism, we get
$$(\Delta_\hbar \otimes \id)(R') {R'}^{1,2}- (\id \otimes \Delta_\hbar)(R'){R'}^{2,3}=0.$$
Let us now define $R={R'}^{-1}$. We have a solution of Theorem \ref{theoprinc}.
\hfill \qed \medskip

\medskip

\section{Concluding remarks}

\medskip

Let us consider the $L_\infty$-quasi-isomorphism of Theorem \ref{theo:map}.
It induces a quasi-isomorphism of complexes:
\begin{theorem}
Let $(\g,\delta_\hbar)$ be a Lie bialgebra. Let $(\Uhg,\Delta_\hbar)$ be
the associated Etingof-Kazhdan quantization. 
The restriction $\varphi$ of the $L_\infty$-quasi-isomorphism of Theorem \ref{theo:map}
defines a quasi-isomorphism between the
coHochschild complex $(T(\Uhg),\dcH)$ (with differential associated to
the coproduct $\Delta_\hbar)$ and $(C(\g),\delta_\hbar)$ the
exterior product over the Lie algebra $\g$ with differential given by
the Lie cobracket.
\end{theorem}
The theorem generalises well known theorem for $\delta_\hbar=0$ (see \cite{Dr:QH},
Proposition 2.2).

\smallskip

Let us discuss some possible generalisations:

\begin{remark}
The main tool used in that paper is Etingof-Kazhdan quantization/dequantiza\-tion
theorem. Suppose now that one can prove an analogue of this theorem for
Lie bialgebroids and Hopf algebroid (which is still a conjecture).
Using the same framework, one would get a 
$L_\infty$-quasi-isomorphism (or even $G_\infty$-quasi-isomorphism,
c.f. Remark \ref{greg}) between the exterior power of a any Lie bialgebroid
and the tensor product of the associated quantized Hopf algebroid.
In the case the Lie algebroid is the algebroid of tensor fields over a manifold $M$,
this would give directly a global formality theorem between tensor fields and
multidifferential operators.
\end{remark}

\begin{remark}
To answer completely to ``Drinfeld last unsolved problem''
(\cite{Dr3}), one should also check the last two conditions of Definition \ref{defi:quanti}: $R^{-1}=R^{2,1}$ and
$R\Delta_\hbar(-)R^{-1}=\Delta_\hbar^{2,1}(-)$.
Those properties are maybe not satisfied for every quantization $\Uhg$ of $(\g,\delta_\hbar)$
\end{remark}

\end{document}